\documentclass[a4paper,12pt]{amsart}
\usepackage{amsmath,amssymb,amscd,graphicx,color}
\newcommand{\R}{\mathbb R}
\newcommand{\C}{\mathbb C}
\newcommand{\Z}{\mathbb Z}
\newcommand{\rpp}{\mathbb R P^2}
\newcommand{\cpp}{\mathbb C P^2}
\newcommand{\cpl}{\mathbb C P^1}
\newcommand{\rpl}{\mathbb R P^1}
\newcommand{\ra}{\mathbb R A}
\newcommand{\ca}{\mathbb C A}
\newcommand{\p}{\partial} 
\renewcommand{\:}{\,{:}\,}
\newcommand{\conj}{\operatorname{conj}}
\newcommand{\inc}{\operatorname{in}}
\newcommand{\Int}{\operatorname{Int}}
\newcommand{\sminus}{\smallsetminus}
\newcommand{\sfit}{\sffamily\slshape}
\theoremstyle{plain}

\newtheorem{lemma}{Lemma}[]

\title[Whitney Number of Real Algebraic Affine Curve of Type I]{Whitney Number of Closed Real Algebraic Affine Curve of Type I}
\author{Oleg Viro}
\begin{document} 
\begin{abstract} 
For a closed real algebraic plane affine curve dividing its
complexification and equipped with a complex orientation, the 
Whitney number is expressed in terms of behavior of its complexification
at infinity.
\end{abstract}

\maketitle

\section{Introduction}\label{s1}
\subsection{Whitney number}\label{s1.1}
Oriented smooth closed immersed curve $C$ on an oriented affine plane has an
important numerical characteristic, {\it Whitney number, \/}  which 
is called also {\it winding number,\/} and can be defined as the rotation 
number of the velocity vector, as well as the degree of the Gauss map 
$C\to S^1$. It determines the
immersion $C\to R^2$ up to  regular homotopy, i.e. path in the space of
immersions.  

\subsection{Real Algebraic Curves Under Consideration}\label{s1.2}
In this paper we consider a class of plane affine real algebraic 
curves such that the Whitney number is defined naturally for their sets
of real points. Namely, we consider irreducible plane affine real 
algebraic curves $A$ satisfying the following three conditions: 
\begin{enumerate} 
\item the set of real points $\ra$ is compact,
\item any real singular point is a non-degenerate double point with 
real branches,
\item the set of real points $\ra$ is zero homologous modulo 2 in the
set of complex points $\ca\subset\cpp$ of the projective closure of $A$. 
\end{enumerate}
The second condition implies that $\ra$ can be
represented as a smoothly immersed curve.  
The first condition implies that $\ra$ is closed, i.e., $A$ has no
real branches going to infinity. Of course, $A$ has complex branches
approaching infinity. They correspond to the intersection points of
$\ca$ and $\cpl_\infty$, the total number of which taken with multiplicities
equals the degree of $A$.

\subsection{Complex Orientations}\label{s1.3}
The third condition is needed to equip $\ra$ 
with a natural orientation. Since $\ra$ is zero homologous in $\ca$
modulo 2, it
bounds a 2-chain : $\ra=\p\ca_+\subset\ca$ modulo 2. Chain $\ca_+$ inherits an
orientation from $\ca$ (which has a natural orientation as a complex
curve), and induces an orientation on $\ra$. Curve $\ra$ bounds also chain 
$\ca_-=\conj\ca_+$, where $\conj:\cpp\to\cpp$ is the complex conjugation
involution, $\conj:(z_0\,{:}\,z_1\:z_2)\mapsto(\bar z_0\,{:}\,\bar
z_1\:\bar z_2)$ . The same construction but using
$\ca_-$ instead of $\ca_+$, defines the opposite orientation on $\ra$. Both
orientations are called {\it complex.\/} If $\ra$ has more than one
connected components, it has other, non-complex orientations. There are
only two chains, $\ca_+$ and $\ca_-$, embedded in $\ca$ and bounded by 
$\ra$. A choice of a complex orientation is equivalent to the choice 
between $\ca_+$ and $\ca_-$.

A real algebraic curve $A$ with $\ra$ zero homologous in $\ca$ is said
to be of {\sfit type I\/}. This definition is due to Felix Klein. 
Any real rational curve with infinite $\ra$ is of type I. More about
curves of type I and complex orientation can be found in \cite{Rokh},
\cite{Viro} and \cite{Zv}.     

For a curve satisfying the conditions of Section \ref{s1.2} and equipped 
with a complex orientation we give an 
interpretation of the Whitney number $w(\ra)$ in terms of behavior of 
a half $\ca_+$ of its complexification $\ca$ at infinity.

\subsection{Line at Infinity}\label{s1.4}
Consider the complex line at infinity $\cpl_\infty=\cpp\sminus\C^2$.
Topologically, this is a 2-sphere. The set $\rpl_\infty=\rpp\sminus\R^2$ 
of its real points is a circle dividing it into two hemi-spheres. 
These hemi-spheres equipped with the orientations inherited from 
$\cpl_\infty$ induce on $\rpl_\infty$ two orientations opposite 
to each other. Denote the
hemi-sphere which defines the positive (counter-clockwise) orientation
on $\rpl_\infty$ by $\cpl_{\infty+}$, and the other hemi-sphere by
$\cpl_{\infty-}$.

\subsection{Main Result}\label{1.4}
{\it 
Let $A$ be a plane affine real algebraic curve satisfying the conditions
of Section \ref{s1.2} and equipped with the complex
orientation defined by $\ca_+$.
Then
\begin{equation}\label{eq1}
w(\ra)=\ca_+\circ\cpl_{\infty+}-\ca_+\circ\cpl_{\infty-}.
\end{equation}
}

In the right hand side of \eqref{eq1}, $\circ$ means intersection
number. 2-Chains $\ca_+$ and $\cpl_{\infty\pm}$ are compact domains of 
complex curves in $\cpp$. Their boundaries $\p\ca_+$ and
$\p\cpl_{\infty+}=\p\cpl_{\infty-}$ do not
meet, there are finitely many points in $\ca_+\cap\cpl_{\infty-}$ and
$\ca_+\cap\cpl_{\infty+}$. Therefore the intersection numbers can be 
defined as the sums of intersection multiplicities over the intersection 
points.

\subsection{Reformulation via Asymptotes}\label{s1.5}
If $\ca$ is transversal to $\cpl_\infty$ then each point of 
$\ca\cap\cpl_\infty$ corresponds to an asymptote of the affine part of 
$\ca$. If the real affine part of this curve is closed, all the
asymptotes are imaginary. Affine imaginary lines which do not meet
$\rpl_\infty$ are divided into those which meet $\cpl_{\infty+}$ and
those which meet $\cpl_{\infty-}$. Theorem \ref{s1.4} claims that $w(\ra)$ 
equals the difference between the number of the asymptotes of $\ca_+$ 
of these two sorts.

\subsection{Sketch of Proof and Organization of Paper}\label{s1.6}
To prove the Main Result, we choose a generic real point on
$\rpl_{\infty}$ and rotate around it oriented real line $L$ counting 
changes of $\ca_+\circ\C L_+-\ca_+\circ\C L_-$. This quantity changes
only when $\R L$ kisses $\R A$. The total change can be identified with
$-2w(\ra)$ calculated as degree of the Gauss map. On the other hand, at
the beginning of rotation, $L$ coincides with $P^1_\infty$, and at the
end, with the same line, but with the opposite orientation.
Therefore, the total change of the quantity is
$-2(\ca_+\circ\cpl_{\infty+}-\ca_+\circ\cpl_{\infty-})$.

For an expert, this sketch would suffice. To make it more formal, we
need to clarify what intersection numbers are to be considered when the
real part of the rotating line $L$ would intersect $\ra$. We have to
exclude intersection points in the real domain, that is on the
boundaries of both 2-chains. To make this in the framework of algebraic
topology, we make a spherical blow up of $\cpp$ along $\rpp$. 

Section \ref{s2} is devoted to blow ups of this kind. In Section \ref{s3} 
the changes of the intersection numbers are calculated, and, in Section 
\ref{s4}, the proof of Main Result is completed. 

\section{Digression on Blowing up of Real Point Set}\label{s2}

\subsection{Blow up a Submanifold}\label{s2.0} For any smooth 
submanifold $Y$ without boundary of a manifold $X$ one
can blow up $X$ along $Y$ in two ways: replacing each point  $y\in Y$ 
either by the projectivization of $T_yX/T_yY$ 
(i.e., the space of real
one-dimensional vector subspaces of $T_yX/T_yY$),
or by the spherization of $T_yX/T_yY$ (i.e., the space of {\it oriented\/}
one-dimensional vector subspaces of $T_yX/T_yY$).  The first kind of
blow up is said to be {\it projective,\/} the second one, {\it
spherical.\/}  A projective blow
up gives a manifold without boundary, while a spherical one gives a
manifold with boundary, which is obtained from $Y$.
If $Y$ is of codimension one in $X$, the projective blow up does not
change $X$, and the spherical blow up cuts $X$ along $Y$, i.e.,
replaces $Y$ by its double covering. 

The set of real points $\ra$ of a non-singular real algebraic variety $A$ is
a smooth submanifold of middle dimension without boundary of the set $\ca$ 
of complex points of $A$. Thus the blow ups outlined above can be
made in this situation.  The specific of situation provides
possibilities for different descriptions of the construction.

Multiplication by $\sqrt{-1}$ defines an isomorphism between $T_y\ra$
and $T_y\ca/T_y\ra$. On the other hand, the projectivization of $T_y\ra$
can be identified with the set of complex one-dimensional subspaces of
$T_y\ca$ invariant under the complex conjugation involution
$T_y\ca\to T_y\ca$. Therefore projective blow up of $\ca$ along $\ra$
can be considered as replacement of each point of $\ra$ by the set of
all real tangent lines of $\ra$ at the point. 

Two orientations of a real line are induced on it as on the boundary of 
the two halves of its complexifications. Therefore, the spherization of 
$T_y\ra$ can be identified with the set of halves of complexifications
of the real lines. The spherical blow up of $\ca$ along $\ra$ can be
identified with replacement of each point of $\ra$ with the set of
halves of complexifications of all real tangent lines of $\ra$ at the
point.

\subsection{Spherical Blow up of Real Projective Space in Complex Projective 
Space}\label{s2.1}

The set of real oriented lines in $n$-dimensional projective space is
naturally identified with the oriented Grassmann variety
$G^+_{2,n-1}(\R)$. Each point $x\in G^+_{2,n-1}(\R)$ is an oriented
2-dimensional vector subspace of $\R^{n+1}$. Its projectivization $Px$ is a
line in $\R P^{n}$ inheriting orientation from $x$.  
Take the set 
$$\Upsilon P^n=\{(l,p)\in G^+_{2,n-1}(\R)\times\C P^n\mid p\in\C Pl_+ \}$$
where $\C Pl_+$ is a hemisphere in the set of complex points of $Pl$ 
such that the orientation of $Pl$ is induced on $Pl$ as on boundary of
$\C Pl_+$ equipped with its complex orientation.
The natural projection 
$$\upsilon:\Upsilon P^n\to\C P^n:(l,p)\mapsto p$$
is bijective over the set of imaginary points because through any
imaginary point one can draw a unique real line (the one that 
is determined by the point and its image under the complex conjugation).
The set of all real oriented lines passing through $p\in\R P^n$  is
homeomorphic to sphere $S^{n-1}$.

Thus $\Upsilon P^n$ can be considered as $\C P^n$  blown up
along $\R P^n$. This is an oriented $2n$-dimensional manifold with boundary. 
The interior of $\Upsilon P^n$ is mapped by $\upsilon$ diffeomorphically 
onto $\C P^n\sminus\R P^n$. The boundary of $\Upsilon P^n$ is mapped by 
$\upsilon$ onto $\R P^n$. The map $\p \Upsilon P^n\to\R P^n$ is a fibration 
with fiber $S^{n-1}$ equivalent to the fibration of unit tangent vectors 
of $\R P^n$, spherization of the tangent bundle of $\R P^n$.

\subsection{Non-Singular Real Projective Variety}\label{s2.2}
The construction of the preceding section is extended naturally to any 
non-singular real algebraic projective variety: for such a variety 
$A\subset P^n$ put
$$\Upsilon A=\{(l,p)\in\Upsilon P^n\mid p\in\C A\text{, and, if }p\in\R
A, \text{ then } l\subset T_p\R A\}.$$

If $A$ is a non-singular real projective curve, $\Upsilon A$ can be 
obtained from $\ca$ by cutting along $\ra$. Recall that cutting of a
surface along a curve two-sidedly embedded into the surface is a 
replacement of the curve by  two copies of it. 

\subsection{Blow up of Real Part in a Singular Curve}\label{s2.3} We need 
the construction of the preceding section only in the case of projective 
plane and a plane curve. However, the curve is not assumed to 
be non-singular. 

To encompass this more general situation, consider, 
for a real plane projective curve $A$, normalization $\nu:\bar A\to A$. 
The set $\C\bar A$ of complex points of $\bar A$ is a compact 
non-singular complex algebraic curve. The restriction of $\conj$ to $\ca$ 
lifts to an anti-holomorphic involution. We will denote it by $c$.
Observe that $\nu^{-1}(\ra)\supset fix(c)$, but there may be points of
$\C\bar A$ which are mapped by $\nu$ to real points of $A$ without being 
fixed under $c$. For example, the preimage under $\nu$ of any isolated 
point of $\ra$ consists of such points.  

For each point $x\in\C\bar A$, the germ of composition 
$$\begin{CD} \bar A @>{\nu}>> A@>{\inc}>>\cpp \end{CD}$$
has a well-defined osculating line. Denote it by $O(x)$.
It passes through $\nu(x)$.
If $x\in fix(c)$ then $O(x)$ is real. It may happen that $x\not\in
fix(c)$, but $O(x)$ is real. 

Denote by $\Upsilon A$ the subset of $G^+_{2,2}(\R)\times\C\bar A $
consisting of pairs $(l,x)$ such that \begin{itemize} 
\item  $Pl=O(x)$, if  $x\in fix(c)$;
\item otherwise just  $\nu(x)\in\C Pl_+$
\end{itemize}

There is a natural map $\Upsilon A\to \C A:(l,x)\mapsto \nu(x)$. On the
preimage of the set of non-singular imaginary points of $A$ it is
bijective, on the preimage of the set of non-singular real points it is
2-1 map.  

\section{Intersection of Complex Halves of Curve and Lines}\label{s3}
\subsection{When Intersection Is Stable}\label{s3.1}
Let $A$ be a plane projective real algebraic curve such that its set of 
real points $\ra$ is zero-homologous modulo 2 in $\ca$. Let
$\ca_+\subset\ca$ be a 2-chain with $\p\ca+=\ra$. Let $L^t$, $t\in\R$ be
a continuous family of real projective lines. Suppose their real point
sets, $\R L^t$, are coherently oriented and the orientations are induced
by the complex orientation of hemispheres $\C L^t_+\subset\C L^t$.

\begin{lemma}\label{L1} 
The number of common imaginary points of $\ca_+$ and $\C L^t_+$
counted with multiplicities, as a function of $t$, is locally
constant at all but finite number of values of $t$, for which $\R L^t$
is either tangent to $\ra$ or passes through singular points of $\ra$. 
\end{lemma}

\begin{proof}Blow up the sets of real points. This gives rise to
oriented compact 2-chains $\Upsilon A$ and $\Upsilon L^t$ in $\Upsilon
P^2$ with $\Upsilon A\cap\p\Upsilon P^2=\p\Upsilon A$ and 
$\Upsilon L^t\cap\p\Upsilon P^2=\p\Upsilon L^t$. 

Each of these 2-chains is covered with two ones, obtained from $\ca_+$ 
and $\ca_-$ in the case of $\Upsilon A$ and 
$\C L^t_+$, $\C L^t_-$ in the case of $\Upsilon L^t$. 
Denote the 2-chain coming from $\ca_+$ and $\C L^t_+$ by 
$\Upsilon A_+$ and $\Upsilon L^t_+$, respectively.    
As above, $\Upsilon A_+\cap\p\Upsilon P^2=\p\Upsilon A_+$ and 
$\Upsilon L^t_+\cap\p\Upsilon P^2=\p\Upsilon L^t_+$

The intersection of
$\Int\Upsilon A_+$ and $\Int\Upsilon L^t_+$ counted with multiplicities
equals the number of common imaginary points of $\ca_+$ and $\C L^t_+$
counted with multiplicities.

If $\R L^t$ does not pass through singular points of $\ra$ and is 
transversal to $\ra$ then $\p\Upsilon L^t_+$ and $\p\Upsilon A_+$ are 
disjoint. Let $U$ be a regular neighborhood of $\p\Upsilon A_+$ in 
$\p\Upsilon P^2$ disjoint from $\p\Upsilon L^t_+$, and $V$ be the 
closure of the complement of $U$ in $\p\Upsilon P^2$. Surfaces 
$\Upsilon A_+$ and $\Upsilon L^t_+$ realize homology classes belonging to
$H_2(\Upsilon P^2, U)$ and $ H_2(\Upsilon P^2, V)$, respectively.
There is intersection pairing 
$$H_2(\Upsilon P^2, U)\times H_2(\Upsilon P^2, V)\to\Z.$$ 
Therefore  $\Int\Upsilon A_+\circ\Int\Upsilon L^t_+=\Upsilon
A_+\circ\Upsilon L^t_+$ can be considered as its value on the homology
classes realized by $\Upsilon A_+$ and $\Upsilon L^t_+$.
The homology class realized by $\Upsilon L^t_+$ does not change under a
sufficiently small change of $t$.
\end{proof}

\begin{lemma}\label{L2}
If, under hypothesis of Lemma \ref{L1}, $A$ satisfies the conditions
listed in Section \ref{s1.2} then the number of common imaginary 
points of $\ca_+$ and $\C L^t_+$ counted with multiplicities jumps only
at $t$ for which $\R L^t$ is tangent to a branch of $\ra$.
\end{lemma}

\begin{proof} According to condition (2) of Section \ref{s1.2}, any
real singular point of $A$ is a non-degenerate double point with real
branches. In the blow up of real part, singularities of this sort are 
resolved.
\end{proof}

\subsection{When Intersection Jumps}\label{s3.2}
Under the assumptions above, assume that
$\R L^{t_0}$ is tangent quadratically to $\ra$ at a (real) point 
$p$ and transversal to $\ca$ at all other intersection points.
Let $t_+$ and $t_-$ be real numbers close to $t_0$ such that 
\begin{itemize}
\item for all $t$
between each of them and $t_0$ line $\R L_t$ is transversal $\ra$, 
\item the number of intersection points of $\R L^{t_+}$ and $\ra$ is
greater by 2 than the number of intersection points of $\R L^{t_-}$ 
and $\ra$.   
\end{itemize} 

\begin{lemma}\label{L3}
If the complex orientations of  $A$ and $L^{t_0}$ at $p$ coincide then
 $$\Int\ca_+\circ\Int\C L^{t_-}_+-\Int\ca_+\circ\Int\C L^{t_+}_+=1 $$
If the complex orientations of $A$ and $L^{t_0}$ are opposite at $p$ then 
   $$\Int\ca_+\circ\Int\C L^{t_-}_+-\Int\ca_+\circ\Int\C L^{t_+}_+=0 $$  
\end{lemma}

This a reformulation of a well-known theorem by Fiedler \cite{F}. 
A sketch of the original proof is presented below. 
The problem under consideration is local: near the point of
tangency the curve and family of lines are standard up to a local
diffeomorphism, which extends locally to complex domain. Therefore it
suffices to prove the lemma for any curve and a family of lines of this 
local diffeomorphic type. For example, we can assume that $A$ is a 
circle and lines $L^t$ are parallel to each other. Intersection 
$\ca\cap\C L^{t_-}$ consists of two points, complex conjugate to each
other. Hemi-sphere $\ca_+$ meets one of the hemi-spheres $\C L^{t_-}_+$ 
and $\C L^{t_-}_-$ and does not meet the other one.

A small perturbation of $A\cup L^{t_-}$ gives a non-singular cubic 
curve with two-component real point set. As an M-curve, it is of type I.
Its set of complex points is obtained by a small perturbation from the 
union of two spheres (which are the sets of complex points of $A$ and 
$L^{t_-}$) meeting each other at two points. The perturbation replaces
a small disk neighborhoods of an intersection point in the spheres with a 
tube connecting the boundary circles of the disks. Therefore, $\ca_+$
meets the half of $\C L^{t_-}$ which gives together with it a half of
the cubic curve.  

In the case, when the complex orientations of $\ra$ and $\R L^{t_0}$  
at the point of tangency coincide, the orientation of the cubic curve
obtained from the complex orientations of $\ra$ and $\R L^{t_-}$
coincides with the complex orientation of the cubic curve. See, for
example, \cite{Rokh}. Hence, in this case
$\ca_+$ meets $\C L^{t_-}_+$.\qed

\section{Proof of the Main Result}\label{s4}

Choose a direction in which $\ra$ has \begin{itemize}
\item no inflection tangent line,
\item no double tangent line,
\item no line tangent to a real branch at a singular point. 
\end{itemize}
Let $o$ be the point on $\rpl_\infty$ which is the common point of all 
lines of the chosen direction.

Consider the pencil $L^t$ of all real lines passing through $o$. Orient them
coherently. Using $L^t$, we evaluate sides of \eqref{eq1}.

Let $t_1$, \dots $t_n$ be the values of $t$ for which $L^t$ are tangent
to $\ra$. They are divided into four classes, (min, $+$), (min, $-$),
(max, $+$), (max, $-$),
according  to the change of the number
of points in $\ra\cap\R L^t$ when $t$ passes $t_i$ and behavior
of orientations of $\ra$ and $\R L^t_i$, see Figure \ref{f1}.

\begin{figure}
\centerline{\input{f1.pstex_t}}
\caption{}
\label{f1}
\end{figure}
 
Obviously, 
$$w(\ra)=\#(min,+)-\#(max,+)=\#(max,-)-\#(min,-).$$

Consider now the change of $\ca_+\circ\C
L^t_+-\ca_+\circ\C L^t_-$ when $t$ passes the critical value.
By Lemma \ref{L3}, in the case (min,$+$)
\begin{multline} 
\Delta(\ca_+\circ\C L^t_+)=-1, \quad \Delta(\ca_+\circ\C L^t_-)=0,\\
\quad \Delta(\ca_+\circ\C L^t_+-\ca_+\circ\C L^t_-)=-1;
\end{multline}
 in the case (min,$-$)
\begin{multline} 
\Delta(\ca_+\circ\C L^t_+)=0, \quad \Delta(\ca_+\circ\C L^t_-)=-1,\\
\quad \Delta(\ca_+\circ\C L^t_+-\ca_+\circ\C L^t_-)=1;
\end{multline}   
 in the case (max,$+$)
\begin{multline} 
\Delta(\ca_+\circ\C L^t_+)=1, \quad \Delta(\ca_+\circ\C L^t_-)=0,\\
\quad \Delta(\ca_+\circ\C L^t_+-\ca_+\circ\C L^t_-)=1
\end{multline}      
 in the case (max,$-$)
\begin{multline} 
\Delta(\ca_+\circ\C L^t_+)=0, \quad \Delta(\ca_+\circ\C L^t_-)=1,\\
\quad \Delta(\ca_+\circ\C L^t_+-\ca_+\circ\C L^t_-)=-1
\end{multline}      

Summing up over all critical lines, we obtain that the total change
of $\ca_+\circ\C L^t_+-\ca_+\circ\C L^t_-$ is equal to
$$-\#(min,+)+\#(min,-)+\#(max,+)-\#(max,-)=-2w(\ra)$$

On the other hand, the pencil $\R L^t$ starts with $\R P^1_\infty$ and
ends up with the same $\R P^1_\infty$, but reverses the orientation.
Therefore it exchanges $\cpl_+$ and $\cpl_-$. Therefore the total change
of $\ca_+\circ\C L^t_+-\ca_+\circ\C L^t_-$ should be equal to
$-2(\ca_+\circ\cpl_{\infty+}-\ca_+\circ\cpl_{\infty-})$.
\qed


\begin{thebibliography}{999}
\bibitem{F} T.~Fiedler, {\em Pencils of lines and topology of real 
algebraic curves}, Izv.  Akad. Nauk, Ser. Mat. {\bf 46} (1982), 853--863 
(Russian), English translation: 
 Math. USSR-Izv. {\bf 21} (1983), 161--170.

\bibitem{Rokh}
V.~A. Rokhlin, {\em Complex topological characteristics of real algebraic
  curves}, Uspekhi Mat. Nauk {\bf 33} (1978), 77--89 (Russian).

\bibitem{Viro}
O.~Ya. Viro, {\em Progress in the topology of real algebraic varieties over the
  last six years}, Uspekhi Mat. Nauk {\bf 41} (1986), 45--67 (Russian), English
  transl., Russian Math. Surveys {\bf41:3} (1986), 55--82.

\bibitem{Zv} V.I.Zvonilov, {\em Complex orientations of real algebraic 
curves with singularities\/}, Soviet Math. Dokl. {\bf 27} (1983), 14--17.

\end{thebibliography}
\end{document}